\documentclass[12pt]{amsart}

\usepackage{amssymb}

\usepackage{enumitem}

\usepackage{hyperref,color}

\makeatletter
\@namedef{subjclassname@2010}{%
  \textup{2010} Mathematics Subject Classification}
\makeatother

\usepackage[T1]{fontenc}

\newtheorem{theorem}{Theorem}[section]

\newtheorem{lemma}[theorem]{Lemma}
\newtheorem{proposition}[theorem]{Proposition}

\theoremstyle{definition}
\newtheorem{definition}[theorem]{Definition}

\numberwithin{equation}{section}

\newcommand{\ep}{\epsilon}
\newcommand{\mbb}{\mathbb}

\newcommand{\ov}{\overline}

\newcommand{\pa}{\partial}
\newcommand{\mf}{\mathbf}

\newcommand{\al}{\alpha}

\newcommand{\z}{\zeta}
\newcommand{\la}{\lambda}
\newcommand{\ti}{\tilde}

\newcommand{\sm}{\setminus}
\newcommand{\ra}{\rightarrow}
\renewcommand{\Re}{\operatorname{Re}}

\hypersetup{
    colorlinks,
    citecolor=blue,
    filecolor=black,
    linkcolor=red,
    urlcolor=black
}

\begin{document}
\title{Further remarks on the higher dimensional Suita conjecture}
\author{G.P. Balakumar, Diganta Borah, Prachi Mahajan and Kaushal Verma}

\address{G. P. Balakumar: Department of Mathematics, Indian Institute of Technology Palakkad, 678557, India}
\email{gpbalakumar@gmail.com}

\address{Diganta Borah: Indian Institute of Science Education and Research, Pune  411008, India}
\email{dborah@iiserpune.ac.in}

\address{Prachi Mahajan: Department of Mathematics, Indian Institute of Technology Bombay, Powai, Mumbai 400076, India}
\email{prachi@math.iitb.ac.in}

\address{Kaushal Verma: Department of Mathematics, Indian Institute of Science, Bangalore 560 012, India}
\email{kverma@iisc.ac.in}


\begin{abstract}
For  a domain $D \subset \mbb C^n$, $n \ge 2$, let $F^k_D(z)=K_D(z)\la\big(I^k_D(z)\big)$, where $K_D(z)$ is the Bergman kernel of $D$ along the diagonal and $\la\big(I^k_D(z)\big)$ is the Lebesgue measure of the Kobayashi indicatrix at the point $z$. This biholomorphic invariant was introduced by 
B\l ocki and in this note, we study its limiting boundary behaviour on two classes of domains namely, $h$-extendible and strongly pseudoconvex polyhedral domains.
\end{abstract}

\maketitle


\section{Introduction}

\noindent We continue the study of $F^k_D$, a biholomorphic invariant that was defined by B\l ocki in his work on Suita's conjecture \cite{B1}. Recall that for a domain $D \subset \mathbb C^n$, 
\[
F^k_D(z) = K_D(z) \la\big(I^k_D(z)\big)
\]
where $K_D(z)$ is the Bergman kernel of $D$ along the diagonal and $\la\big(I^k_D(z)\big)$ is the Lebesgue measure of the Kobayashi indicatrix at $z \in D$. As usual, let $k_D = k_D(z, v)$ be the infinitesimal Kobayashi metric on $D$. 
B\l ocki--Zwonek \cite{BZ} have shown that
\[
1 \le F^k_D(z) \le C^n
\] 
where $C = 4, 16$ accordingly as $D$ is convex or $\mathbb C$-convex respectively. Furthermore, their work also contains a detailed discussion of this invariant on convex egg domains in $\mathbb C^2$. These results were supplemented in \cite{BBMV} 
wherein this invariant was considered on strongly pseudoconvex domains in 
$\mathbb C^n$ and a few other observations were made about its {\it boundary behaviour} on egg domains in $\mathbb C^2$. In particular, even on the smoothly bounded convex eggs of the form
\[
E_{2 \mu} = \big\{(z, w) \in \mathbb C^2: \vert z \vert^2 + \vert w \vert^{2 \mu} < 1 \big\}
\]
for integers $\mu>1$, $F^k_{E_{2 \mu}}$ does not admit a limit at any of the weakly pseudoconvex points of $\partial E_{2 \mu}$. In fact, the full range of all possible values of  $F^k_{E_{2 \mu}}$ at points of $E_{2 \mu}$ show up as possible limits near any of the weakly pseudoconvex points on $\partial E_{2 \mu}$. By the well known work of Lempert, all invariant metrics on bounded convex domains $D$ coincide; so in particular for any
invariant metric $\tau$, $F_D^k \equiv F_D^\tau$, the analogously defined invariant function
associated to $\tau$. While such an identity need not hold on strongly pseudoconvex domains in general, it was shown in \cite{BBMV} that on any smoothly bounded strongly pseudoconvex domain $D$, the boundary limits of $F_D^\tau$ exist and give rise to the same value: $F^\tau_D(z) \rightarrow 1$ as $z$ approaches $\partial D$, the boundary of $D$. A different approach to this result has been suggested recently in \cite{BZr}, wherein the focus was the invariant metric $\tau=a$ of Azukawa. In this article, 
B\l ocki -- Zwonek have also raised questions about the boundary behaviour of $F_D^a$ both for bounded convex domains as well as for smoothly bounded pseudoconvex domains. While the aforementioned (convex) egg domains settle the non-existence of boundary limits of $F_D^a$ at non-strongly pseudoconvex boundary points in general, it is possible to make certain definite statements about the possible limiting boundary values. 

\medskip

The purpose of this note is to record some general properties of $F^k_D$ and to compute its possible limiting boundary values on $h$-extendible and strongly pseudoconvex polyhedral domains.

\section{Some Observations}

\noindent (i) \textit{Removable Singularities}: For a bounded domain $G \subset \mathbb C^n$ and a subvariety $V \subset G$ of codimension at most $2$, it is known that $k_G(z, v) = k_{G \sm V}(z, v)$ for $(z, v) \in (G \sm V) \times \mathbb C^n$. Further, the Bergman kernels along the diagonal of $G$ and $G \sm V$ are equal since $V$ is a removable singularity for $L^2$-holomorphic functions. Hence $F^k_G = F^k_{G \sm V}$.

\medskip 

If $V$ has codimension $1$, this is no longer the case despite the fact that $V$ is still removable for $L^2$-holomorphic functions. In general, $k_G \le k_{G \setminus V}$ as can be seen by taking $G$ to be the unit disc $\mathbb{D} \subset \mbb C$ and $V = \{0\}$. Thus, the most that can be said in general is that $F^k_G \le F^k_{G \sm V}$. Examples in higher dimensions can be constructed by observing that  $F^k_D$ is multiplicative as a function of the domain $D$, that is, 
\[
F^k_{D \times G}\big((p, q)\big) = F^k_D(p)F^k_G(q)
\]
for $D \subset \mbb C^n, G \subset \mbb C^m$, and hence $F^k_{D \times \mbb D} \le  F^k_{D \times \left(\mbb D\setminus \{0\}\right)}$. Finally, using the
multiplicative property, this invariant for the Hartogs' triangle $\Omega \subset \mbb C^2$ (which is biholomorphic to $\mbb D \times \left(\mbb D \setminus \{0\}\right))$ can be computed as
\[
F^k_{\Omega}(z,w)=F^k_{\mbb D}\left( \frac{w}{z} \right) F^k_{\mbb D \setminus \{0\}}(z)  = 4 \left(\frac{\vert z \vert \log \vert z \vert}{1-\vert z \vert^2}\right)^2.
\]
In particular, $F^k_{\Omega}(z,w) \to 0$ as $(z,w)$ approaches the origin from within $\Omega$.

\medskip

\noindent (ii) \textit{Regularity of $F_G^k$}: Let $G \subset \mbb C^n$ be an arbitrary domain. Then $F_G^k(z)$ is always lower semicontinuous; it is continuous when the Kobayashi metric 
$k_G(z, v)$ is continuous on $G \times \mbb C^n$ and non-degenerate in the sense that $k_G(z, v) > 0$ for all $v \not= 0$.

\medskip

It suffices to show lower semicontinuity of the function $z \mapsto \lambda(I^k_G(z))$. For this, fix a $z^0 \in G$ and let $z^{\nu}$ be a sequence in $G$ converging to $z^0$. We claim that
\begin{equation}\label{I-z0-znu}
I^k_G(z^0) \subset \liminf_{\nu \to \infty} I^k_G(z^{\nu})= \cup_{\mu= 1}^{\infty} \cap_{\nu= \mu}^{\infty} I^k_G(z^{\nu}).
\end{equation}
Indeed, let $v^0 \in I^k_G(z_0)$. Then $\ep=1 - k_G(z^0,v^0)>0$ and by the upper semicontinuity of $k_G(\cdot,v^0)$, there exists $N \in \mathbb{N}$ (depending possibly on both $z^0$ and $v^0$) such that for all $\nu \geq N$ we have
\[
k_G(z^{\nu},v^0) \leq k_G(z^0,v^0) + \epsilon/2 = 1-\ep/2<1.
\]
This implies that $v^0 \in I^k_G(z^{\nu})$ for all $\nu \geq N$ proving our claim \eqref{I-z0-znu}. Now by Fatou's lemma for measurable sets,
\[
\lambda\big(I^k_G(z^0)\big) \leq \lambda\big(\liminf I^k_G(z^\nu)\big) \leq \liminf_{\nu \to \infty} \lambda\big(I^k_G(z^\nu)\big)
\]
which establishes the lower semicontinuity of $z \mapsto \lambda\big(I^k_G(z)\big)$. 

Now, we restrict attention to those domains $G$ for which $k_G$ is continuous and non-degenerate. To show the continuity of $F^k_G$, it suffices to show that $z \mapsto \lambda\big(I^k_G(z)\big)$ is upper semicontinuous. For this, pick a $z^0 \in G$ and let $z^{\nu}$ be a sequence in $D$ converging to $z^0$. Let $\ep>0$. By continuity of $k_G$, there exists $N \in \mf{N}$ such that
\[
k_G(z^{\nu},v) > k_G(z^0,v) - \epsilon
\]
for all $\nu \geq N$ and $v \in S = \partial \mbb B^n$, the standard Euclidean unit sphere. Now, for any nonzero vector $v \in I^k_G(z^{\nu})$ where $\nu \geq N$,
\[
1 > k_G(z^{\nu},v) = \vert v \vert k_G(z^{\nu}, v/\vert v \vert) > \vert v \vert \big( k_G(z^0,v/\vert v \vert) - \epsilon \big) = \vert v \vert \big( k_G(z^0,v)/\vert v \vert - \epsilon \big)
\]
and thus
\[
1 + \epsilon \vert v \vert > k_G(z^0,v).
\]
The continuity of $k_G$ together with its non-degeneracy also implies that there is a positive constant $c$ depending only on $z^0$ such that $k_G(z^{\nu},v) \geq c \vert v \vert$ for all $\nu \geq N$ and $v \in \mathbb{C}^n$, modifying $N$ if necessary. This implies that $I^k_G(z^{\nu}) \subset c^{-1} \mbb B^n$ for all $\nu \geq N$. 

It follows that any vector $v$ picked from any of the indicatrices $I^k_G(z^{\nu})$ for $\nu \geq N$ satisfies
\[
k_G(z^0,v) < 1 + c^{-1}\epsilon
\]
and this means that $I^k_G(z^{\nu}) \subset (1+c^{-1}\epsilon)I^k_G(z^0)$. Therefore, 
\[
\lambda\big(I^k_G(z^{\nu})\big) \leq (1 + c^{-1}\epsilon)^{2n}  \lambda\big(I^k_G(z^0)\big)
\]
for all $\nu \geq N$. As $c$ depends only on $z^0$ but not on $\epsilon$, this implies
\[
\limsup_{\nu \to \infty} \lambda\big(I^k_G(z^{\nu})\big) \leq \lambda\big(I^k_G(z^0)\big)
\]
proving the upper semicontinuity of $z \mapsto \lambda\big(I^k_G(z)\big)$.

\medskip
 
To conclude, note that the non-degeneracy of $k_G$ entails the boundedness of the indicatrices and its continuity implies the following geometric property: every ray emanating from the origin in $T_z(G)$ intersects the (usual topological) boundary $\partial I_G^k(z)$ of the indicatrix, in exactly one point. This then leads to $\partial I_G^k(z)$ being homeomorphic 
to $\partial \mathbb{B}^n$; indeed, as the `graph' of a (uniformly) continuous function on $\partial \mathbb{B}^n$ thereby ensuring that the $2n$-dimensional Lebesgue measure of $\partial I_G^k(z)$ is zero. 

\medskip

\noindent (iii) \textit{A remark on $F_{\mbb G_2}^k$}: For the symmetrized bidisc $\mbb G_2 \subset \mbb C^2$, it is known that $F^k_{\mbb G_2}\big((0, 0)\big) = 4/3$. The explicit description of ${\rm Aut}(\mbb G_2)$ shows that every point in $\mbb G_2$ is equivalent to $(b, 0)$ for $0 \le b < 1$. A rough estimate of $F^k_{\mbb G_2}\big((b, 0)\big)$ can be obtained by using the known expression for the Bergman kernel $K_{\mbb G_2}\big((b, 0)\big)$ and then estimating the volume of $I^k_{\mbb G_2}\big((b, 0)\big)$.

\medskip

Indeed,
\[
K_{\mbb G_2}\big((b, 0)\big) = \frac{2-b^2}{\pi^2(1-b^2)^2}
\]
and if $\mbb G_2 = \pi(\mbb D^2)$ where $\pi(z_1, z_2) = (z_1 + z_2, z_1 z_2)$ is the symmetrizing map, then 
\[
d \pi \big((b, 0)\big) : I^k_{\mbb D^2}\big((b, 0)\big) \rightarrow I^k_{\mbb G_2}\big((b, 0)\big) 
\]
and hence
\[
\lambda\Big(I^k_{\mbb G_2}\big((b, 0)\big)\Big) \ge \Big\vert \det d \pi \big((b, 0)\big) \Big\vert^2 \lambda\Big(I^k_{\mbb D^2}\big((b, 0)\big)\Big).
\]
Since $\big\vert \det d \pi \big((b, 0)\big) \big\vert^2 = b^2$ and $\lambda\big(I^k_{\mbb D^2}\big((b, 0)\big) \big)= \pi^2(1-b^2)^2$, it follows that
\[
F^k_{\mbb G_2}\big((b, 0)\big) \ge (2-b^2)b^2
\]
for $0 \le b < 1$. This lower bound does not yield anything interesting near $b = 0$ for then it is expected that $F^k_{\mbb G_2}\big((0, 0)\big) \approx 4/3$, but it does show that $\liminf_{b \ra 1} F^k_{\mbb G_2}\big((b, 0)\big) \ge 1$. Finally, note that $\pi^{-1}\big((b, 0)\big) = \big\{(b, 0), (0, b)\big\}$ if $b > 0$ and one could work with $d \pi\big((0, b)\big)$ which also maps $I^k_{\mbb D^2}\big((0, b)\big)$ into $I^k_{\mbb G_2}\big((b, 0)\big)$. However, this gives the same lower bound for $F^k_{\mbb G_2}\big((b, 0)\big)$ as above. It would be interesting to see if there is a transformation formula for the  Kobayashi indicatrix under proper holomorphic maps -- this may lead to better estimates for $F^k_{\mbb G_2}\big((b, 0)\big)$. 
 
\medskip

\noindent (iv) \textit{Localization:} It is possible to localize this invariant near peak points as follows:
\begin{proposition} \label{localisation2}
 Let $ G \subset \mathbb{C}^n $ be a pseudoconvex domain and let $ p \in \partial G$ be a local holomorphic peak point. Then for a sufficiently small 
neighbourhood $ U $ of $ p $,  
\[
 \lim_{U \cap G \ni z \rightarrow p} \frac{F^{k}_{U \cap G} (z)}{F^{k}_G(z)} = 1.
\]
\end{proposition}
It should be mentioned that this holds for $F_D^\tau$ where $ \tau=a$ the Azukawa metric as well; the proof is immediate when the
already known localization properties of the Azukawa and the Kobayashi metrics 
(cf. \cite{Nik}, \cite{Nik2}, \cite{Graham}) are combined with that of the Bergman kernel (cf. \cite{Hor}, \cite{Nik}). 

\section{$h$-extendible domains}

\noindent Recall that a  pseudoconvex domain $D \subset \mbb C^{n+1}$ is said to 
be $h$-extendible near a smooth, finite-type point $p \in \partial D$ if the Catlin multitype $(1, m_1, \ldots, m_n)$ of $\pa D$ at $p$ satisfies $m_{n-q+1} = \Delta_q < \infty$ for $1 \le q \le n$, where $\Delta_q$ is the $q$-type of $p$. In this case, there are local coordinates $z = (z_0, z') = (z_0, z_1, \ldots, z_n)$ around $p =0$ and a real-valued, plurisubharmonic, $(1/m_1, 1/m_2, \ldots, 1/m_n)$ weighted homogeneous polynomial $P$ of total weight $1$ with no pluriharmonic terms such that $D$ is defined locally near $p$ by
\begin{equation}
\rho(z) = \Re z_0 + P(z', \overline z') + R(z)
\end{equation}
where $R(z) \lesssim \left(\vert z_0 \vert + \vert z_1 \vert^{m_1} + \ldots + \vert z_n \vert^{m_n} \right)^{\gamma}$  for some $\gamma > 1$. Call
\[
D_{\infty} = \big\{(z_0, z') : \Re z_0 + P(z') < 0 \big\}
\]
the local model for $D$ at $p$. Then $D_{\infty}$ is a taut domain. It is known that $h$-extendability of $D$ at $p$  is equivalent to the existence of a positive, $C^{\infty}$-smooth function $a(z')$ on $\mbb C^n \sm \{0\}$ such that $a$ is weighted homogeneous with the same weights as for $P$ and $P(z') - \epsilon a(z')$ is strictly plurisubharmonic on 
$\mbb C^n \sm \{0\}$ when $0 < \epsilon \le 1$. Note that Levi co-rank one, convex finite type and decoupled finite type domains are all examples of $h$-extendible domains. More details can be found in \cite{BSY} and \cite{Yu}.

\begin{theorem}
Let $D \subset \mbb C^{n+1}$ be a bounded pseudoconvex domain that is $h$-extendible at $p \in \pa D$ with multitype $(1, m_1, \ldots, m_n)$ and whose associated local model is $D_{\infty}$. If $\Gamma$ is a non-tangential cone in $D$ with vertex at $p$, then
\[
\lim_{\Gamma \ni z \ra p} F^k_D(z) = F^k_{D_{\infty}}(b)
\]
where $b = (-1, 0, \ldots, 0) \in D_{\infty}$.
\end{theorem}

\begin{proof}
The boundary behaviour of $K_D(z)$ as $z \ra 0$ within $\Gamma$ is known. Indeed, Theorem 1 in \cite{BSY} shows that
\[
\lim_{\Gamma \ni z \ra 0} K_D(z) \big\vert \rho(z) \big\vert^{\beta} = K_{D_{\infty}}(b)
\]
where $\beta = \sum_{j = 0}^n 2/m_j$.

\medskip

To handle the Kobayashi indicatrices, first fix an $\ep \in (0, 1)$ and let $U_{\ep}$ be a neighbourood of $p = 0$ such that the bumped model
\[
D_{\ep} = \big\{ \Re z_0 + P(z') - \ep a(z') < 0 \big\}
\]
contains $D \cap U_{\ep}$. By \cite{Yu}, $D_{\ep}$ is taut. For $t > 0$, let
\[
\pi_t(z) = (tz_0, t^{1/m_1} z_1, \ldots, t^{1/m_n} z_n)
\]
and note that the scaled domains $D_z = \pi_{1/\vert \rho(z) \vert}(D \cap U_{\ep})$ converge to $D_{\infty}$ in the Hausdorff sense. Also, if $z \ra 0$ within $\Gamma$, the base points $\zeta(z) = \pi_{1/\vert \rho(z) \vert}(z)$ converge to a compact subset of the line $\{\Re z_0 = -1, z' = 0\} \subset D_{\infty}$. This is so since non-tangential convergence implies that $\vert \Re z_0 \vert \approx \vert \rho(z) \vert$. Finally, note that $\pi_t \in {\rm Aut}(D_{\infty})$ for all $t > 0$ and hence $D_z \subset D_{\ep}$ for all $z \in \Gamma$ close to the origin. Theorem 2.1 of \cite{Yu} shows that for all fixed $w \in D_{\infty}$, the Kobayashi metrics $k_{D_z}(w, v) \ra k_{D_{\infty}}(w, v) > 0$ as $z \ra 0$ within $\Gamma$. Moreover, the convergence is uniform on compact subsets of $D_{\infty} \times \mbb C^n$. Hence, the indicatrices $I^k_{D_z}(w)$ converge to $I^k_{D_{\infty}}(w)$ in the Hausdorff sense and $\lambda\big(I^k_{D_z}(w_j)\big) \ra \lambda\big(I^k_{D_{\infty}}(w_0)\big)$ if $w_j \ra w_0 \in D_{\infty}$. In particular, as $z \ra 0$ within $\Gamma$,
\[
\lambda\Big(I^k_{D_z}\big(\zeta(z)\big)\Big) \ra \lambda\big(I^k_{D_{\infty}}(\tilde \zeta)\big)
\]
where $\tilde \zeta$ is a possible limit point of $\zeta(z)$. But as noted above, $\tilde \zeta$ lies on the line $\{\Re z_0 = -1, z' = 0\}$ and since $D_{\infty}$ is invariant under volume preserving translations of the form $z \mapsto z + i \alpha$, $\alpha \in \mbb R$, it follows that 
$\lambda\big(I^k_{D_{\infty}}(\tilde \zeta)\big) =  \lambda\big(I^k_{D_{\infty}}(b)\big)$.

\medskip

To conclude, it remains to note that 
\[
\lambda\Big(I^k_{D_z}\big(\zeta(z)\big)\Big) = \big\vert \rho(z) \big\vert^{-\beta} \lambda\big(I^k_D(z)\big)
\]
by the change of variables formula and that
\[
F^k_D(z) = K_D(z) \big\vert \rho(z) \big\vert^{\beta} \lambda\big(I^k_D(z)\big) \big\vert \rho(z) \big\vert^{-\beta}  \ra K_{D_{\infty}}(b) \lambda\big(I^k_{D_{\infty}}(b)\big) = F^k_{D_{\infty}}(b) 
\]
as $z \ra 0$ within $\Gamma$.
\end{proof} 
It should be noted that the non-tangential condition cannot be dropped as the example of $F^k_{E_{2 \mu}}$ shows. More precisely, for $D= E_{2 \mu}$ and $q$ one of its weakly pseudoconvex points (say $q=(0,1)$), note that 
\[
D_\infty = \{ (z,w) \in \mathbb{C}^2 \; : \; 2 \Re(z) +
 \vert w \vert^{2 \mu} <0 \}.
\]
Indeed, an analogue of the Cayley transform maps $D_\infty$ biholomorphically onto $E_{2 \mu}$ with
$b=(-1,0)$ mapped to the origin where the value of $F^k_{E_{2 \mu}}$ is $1$. So,
\[
\lim_{ \Gamma \ni z \to q} F^k_{E_{2 \mu}} =1
\]
whereas we know from \cite{BBMV} that every value attained by $F^k_{E_{2 \mu}}(z)$ as 
$z$ varies in $E_{2 \mu}$ is also attained as a boundary limiting value at $q$; in particular (as $F^k_{E_{2 \mu}}$ is a constant function only for $\mu = 1$), there are sequences
$\{ p_n\} \subset E_{2 \mu}$ approaching $q$ non-tangentially along which $F^k_{E_{2 \mu}}$ has a limit and, the limiting value differs from $1$ -- for instance, follow any one particular orbit of a point of the 
form $(0,p)$ with $0<p<1$, under the action of the automorphism group ${\rm Aut}(E_{2 \mu})$. But then, it turns out that `highly tangential sequences' again yield boundary limit $1$ -- to record this peculiar feature of the boundary behaviour of $F^k_{E_{2 \mu}}$ at the weakly pseudoconvex points of
$F^k_{E_{2 \mu}}$ a bit more precisely but briefly, let $q_n$ be a sequence of points
in $E_{2 \mu}$ converging to $q$ such that: (i) the inner products of $(q_n -q)/\vert q_n -q \vert$
with the unit inner normal to $\partial E_{2 \mu}$ converge to $0$ and 
(ii) the $q_n$'s belong to mutually distinct orbits of ${\rm Aut}(E_{2 \mu})$. Then,
$\lim_{n \to \infty} F^k_{E_{2 \mu}}(q_n) =1$.\\

\noindent It is possible to obtain global bounds for this invariant on Levi co-rank one domains. This follows from the following wherein $\tau$ denotes any distance decreasing metric or the Bergman metric.

\begin{lemma} \label{Levico1}
Let $\Omega$ be a smoothly bounded pseudoconvex domain in $\mathbb{C}^n$ whose boundary  
$\partial \Omega$ is of (finite type and of) Levi corank at most one at $p \in \partial \Omega$. Then there exist positive constants $c,C$ and a neighborhood $U$ of $p$ such that 
$c \leq F_\Omega^\tau(z) \leq C$ for all $z\in \Omega\cap U$.
\end{lemma}

\begin{proof}
This will follow from the well-known boundary estimates of Catlin and 
Cho for $\tau$ equal to the Carath\'{e}odory, Kobayashi or the Bergman metric. We recall the relevant ideas briefly. Let $r$ be a local defining function for $\pa \Omega$ in a neighbourhood $U$ of $p=0$. By shrinking this neighbouhood if needed, we may assume that the orthogonal projection onto the boundary $\partial D$ is well-defined on $U$ and that the normal vector field, given at any $\zeta \in U$ by
\[
\nu(\zeta) = \big(\pa r /\pa \ov z_1(\zeta), \pa r / \pa \ov z_2(\zeta), \ldots, \pa r / \pa \ov z_n(\zeta) \big)
\]
has no zeros in $U$; this is normal to the hypersurface $\Gamma_\zeta = \big\{ r(z) = r(\zeta) \big\}$.

\medskip

\noindent For each  
$\z \in U$, there is a uniform radius $R > 0$ and an injective holomorphic mapping $\Phi^{\z} : B(\z, R) \ra \mbb C^n$ such that the transformed defining 
function $\rho^{\z} = r^{\z} \circ (\Phi^{\z})^{-1}$ reads
\begin{multline}\label{nrmlfrm}
\rho^{\z}(w) = r(\z) + 2 \Re w_n + \sum_{l = 2}^{2m}P_l(\z; w_1) + \vert{w_2}\vert^2 + \ldots + \vert{w_{n-1}}\vert^2 \\
+ \sum_{\al = 2}^{n - 1} \sum_{ \substack{j + k \le m\\
                                                          j, k > 0}} \Re \Big( \big(b_{jk}^{\al}(\z) w_1^j \ov w_1^k \big) w_{\al} \Big)+ R(\z;w)
\end{multline}
where 
\[
P_l(\zeta; w_1) = \sum_{j + k = l} a^l_{jk}(\zeta) w_1^j \ov w_1^k
\]
are real valued homogeneous polynomials of degree $l$ without harmonic terms and the error function $R(\zeta,w) \to 0$ as $w \to 0$ faster than one of the monomials of weight $1$. Further, the map $\Phi^{\z}$ is 
actually a holomorphic polynomial automorphism of weight one  of the form
\begin{alignat}{3} \label{E45}
\Phi^{\z}(z) = \Big(z_1 - \z_1, G_{\z}(\ti z - \ti \z) - Q_2(z_1 - \z_1), \langle \nu(\zeta), z- \zeta \rangle - Q_1({}'z - {}' \zeta) \Big)
\end{alignat}
where $G_{\zeta} \in GL_{n-2}(\mbb C), \ti z = (z_2, \ldots z_{n-1}), {}'z = (z_1, z_2, \ldots, z_{n-1})$ and $Q_2$ is a vector valued polynomial whose $\al$-th component is a polynomial 
of weight at most $1/2$ of the form
\[
Q^{\al}_2(t) = \sum_{k = 1}^m b_k^{\al}(\z) t^k
\]
for $t \in \mbb C$ and $2 \le \al \le n - 1$. Finally, $Q_1({}' z - {}'\zeta)$ is a polynomial of weight at most $1$ and is of the form $\hat{Q}_1 \big(z_1 -\zeta_1, G_\zeta (\tilde{z} - \tilde{\zeta}) \big)$ with $\hat{Q}_1$ of the form
\[
\hat{Q}_1(t_1, t_2, \ldots, t_{n-1}) = \sum_{k=2}^{2m} a_{k0}(\z)t_1^k - \sum_{\al = 2}^{n-1} \sum_{k=1}^m a_k^{\al}(\z) t_{\al}t_1^k - \sum_{\al = 2}^{n-1}c_{\al}(\z) t^2_{\al}.
\]
Since $G_\zeta$ is just a linear map, $Q_1({}'z - {}' \zeta)$ also has the same form when considered as an element of the algebra of holomorphic polynomials 
$\mathbb{C}[{}'z - {} '\zeta]$, when $\zeta$ is held fixed. The coefficients of all the polynomials, mentioned above, are smooth functions of $\zeta$. Note that $\Phi^{\z}(\z) = 0$ and 
\[
\Phi^{\z}(\z_1, \ldots, \z_{n-1}, \z_n - \ep) = \big(0, \ldots, 0, -\ep \;\pa r/ \pa \ov z_n(\z)\big).\\
\]

\noindent Define for each $\delta > 0$, the special-radius
\begin{alignat}{3} \label{E46}
\tau(\z, \delta) = \min \Big\{ \Big( \delta/ \vert P_l(\zeta, \cdot) \vert \Big)^{1/l}, \; \Big(\delta^{1/2}/B_{l'}(\zeta) \Big)^{1/l'}  \; : \; 2 \le l \le 2m, \; 2 \leq l' \leq m   \Big \}.
\end{alignat}
where
\[
B_{l'}(\zeta) = {\rm max} \Big\{ \big\vert b_{jk}^\alpha (\zeta) \big\vert \; : \; j+k=l', \; 2 \leq \alpha \leq n-1 \Big\}, \; 2 \leq l' \leq m.
\]
Here, the norm of the homogeneous polynomials $P_l(\zeta, \cdot)$ of degree $l$, is taken according to the following convention: for a homogeneous polynomial 
\[
p(v) = \sum\limits_{j+k=l} a_{j,k} v^j \bar{v}^k,
\]
define $\vert p(\cdot) \vert = \max_{\theta \in \mathbb{R}} \vert p(e^{i \theta}) \vert$. It was shown in \cite{Cho2} that the coefficients $b_{jk}^\alpha$'s in the above definition of $\tau(\z,\delta)$ are insignificant and may be ignored, so that 
\[
\tau(\z, \delta) = \min \Big\{ \Big( \delta/ \vert P_l(\zeta, \cdot) \vert \Big)^{1/l}  \;: \; 2 \le l \le 2m \Big\}.
\]
Set
\[
\tau_1(\z, \delta) = \tau(\z, \delta) = \tau, \tau_2(\z, \delta) = \ldots = \tau_{n-1}(\z, \delta) = \delta^{1/2}, \tau_n(\z, \delta) = \delta
\]
and define the dilations
\[
\Delta_\zeta^\delta(z) = \big(z_1/\tau_1(\z, \delta), \ldots, z_n/\tau_n(\z, \delta)\big).
\]
The scaling maps are defined by the composition 
\[
S_\zeta^\delta(z) = \Delta_\zeta^\delta \circ \Phi^\zeta.
\]
and they induce the so-called $M$-metric defined on the one-sided neighbourhood $U \cap D$:
\[
M_D(\zeta, v) = \sum\limits_{k=1}^{n} \Big\vert \big(D \Phi^\zeta(\zeta)v\big)_k \Big\vert \big/\Big\vert \tau_k\big(\zeta, \epsilon(\zeta)\big)\Big\vert =  \Big\vert D\big(S_\zeta^\delta(\zeta)\big) (v) \Big\vert_{l^1}
\]
where $\epsilon(\zeta)>0$ is such that $\tilde{\zeta} = \zeta + (0, \ldots, \epsilon(\zeta)$ lies on $\partial D$. The significance of this metric is that it is uniformly comparable to the Kobayashi metric \cite{TT}, in the sense that 
\begin{equation}\label{Kobest}
k_D(\zeta, v)  \approx M_D(\zeta,v) \approx 
\Big\Vert D\big(S_\zeta^\delta(\zeta)\big) (v) \Big\Vert,
\end{equation}
where $\Vert \cdot \Vert$ denotes any norm on $\mathbb{C}^n$ and the suppressed constants are independent of $v$ and $\zeta$ (depending only on the domain 
$U \cap D$). In particular, taking 
$\Vert \cdot \Vert$ to be the $l^\infty$-norm, we may translate this estimate on the Kobayashi metric into one about its indicatrix and its dilates:
\begin{equation} \label{Ball-Box}
c R(\zeta) \subset I^k_D(\zeta)  
\subset C R(\zeta)
\end{equation}
here $c, C$ are a pair of positive constants independent of $\zeta$ and $R(\zeta)$ is the polydisc centered at the origin of polyradius 
\[
\Big(\tau_1\big(\zeta, \ep(\zeta)\big), \sqrt{\ep(\zeta)}, \ldots, \sqrt{\ep(\zeta)}, 
\ep(\zeta) \Big).
\]
When $D$ is additionally bounded and globally pseudoconvex, it follows from  Theorem 1 of \cite{cho1} that for all $\zeta$ in some tubular neighborhood of $U \cap \partial D$, 
\[
K_D(\z, \overline{\z}) \approx \Big( {\rm Vol} \big( R(\zeta) \big) \Big)^{-1}
\]
wherein the suppressed constants depend only on $D$ and are independent of $\zeta$. Combining this with (\ref{Ball-Box}), finishes the verification of Lemma \ref{Levico1}.
\end{proof}
We now further note that as a consequence of Proposition \ref{localisation2}, the boundedness 
assumption on $\Omega$ in the above lemma can be dropped, provided we restrict to $\tau= a$ or $k$, the 
Azukawa or Kobayashi metrics respectively. Indeed, to see that the global smoothness assumption of the
above lemma can be circumvented when we combine it with Proposition \ref{localisation2} (for $\tau=a,k$) to
drop the boundedness assumption, we recall a technique explained by Bell in the final section of his 
article \cite{Bell}, for completing a small piece of $\partial \Omega$ (in case $\Omega$ is unbounded) 
about $p$ into a smooth pseudoconvex finite type hypersurface so that the resulting {\it smoothly} bounded 
domain $G$ is a (small) subdomain of $\Omega$; the lemma above applies to $G$ and then, Proposition \ref{localisation2} will compare $F_G^\tau$ 
with $F_\Omega^\tau$ to yield the version of the above lemma for unbounded $\Omega$ as desired in case $\tau=a,k$. Further next, the just-mentioned technique of Bell, also enables us to 
drop the global smoothness and boundedness assumption to conclude the following
\begin{theorem}
Let $D \subset \mbb C^{n+1}$ be a pseudoconvex domain whose boundary is smooth and of finite type near $p \in \partial D$. Suppose $p$ is an $h$-extendible point with $D_{\infty}$ being the associated local model. If $\Gamma$ is a non-tangential cone in $D$ with vertex at $p$, then
\[
\lim_{\Gamma \ni z \ra p} F^k_D(z) = F^k_{D_{\infty}}(b)
\]
where $b = (-1, 0, \ldots, 0) \in D_{\infty}$.
\end{theorem}


\section{piecewise smooth strongly pseudoconvex domains}

\begin{definition} \label{ps_defn}
A bounded domain $ D $ in $ \mathbb{C}^n $ is said to be a strongly pseudoconvex polyhedral domain with piecewise smooth boundary if 
there are $ C^2$-smooth real valued functions 
$ \rho_1, \ldots, \rho_k : \mathbb{C}^n \rightarrow \mathbb{R} $, $ k \geq 2 $ such that
\begin{enumerate}
 \item[(i)] $ D = \big\{ z \in \mathbb{C}^n: \rho_1(z) < 0, \ldots, \rho_k(z) < 0 \big\} $,
 \item[(ii)] for $ \{ i_1, \ldots, i_l \} \subset \{1, \ldots, k\} $, the gradient vectors $ \nabla \rho_{i_1}(p), \ldots, \nabla \rho_{i_l}(p) $ are linearly independent 
 over $ \mathbb{C} $ for every point $ p $ such that $ \rho_{i_1}(p) = \ldots = 
 \rho_{i_l}(p) = 0 $, and,
 \item[(iii)] $ \partial D $ is strongly pseudoconvex at every smooth boundary point,
\end{enumerate} 
\end{definition}
where for each $ i = 1, \ldots, k $ and $ z \in \mathbb{C}^n $,
\[
 \nabla \rho_i (z) = 2 \left( \frac{\partial \rho_i}{\partial \bar{z}_1}(z), \ldots, \frac{\partial \rho_i}{\partial \bar{z}_n}(z) \right).
\]
Since the intersection of finitely many domains of holomorphy is a domain of holomorphy, it follows that the polyhedral domain $ D $ as in Definition \ref{ps_defn} is 
pseudoconvex. 

\medskip

Let $ D \subset \mathbb{C}^2 $ be a strongly pseudoconvex polyhedral domain with piecewise smooth boundary as above defined by
\[
 D = \big\{ z \in \mathbb{C}^2: \rho_1(z) < 0, \ldots, \rho_k(z) < 0 \big\}. 
\]
Let $ p^0 \in \partial D $ be a singular boundary point, i.e., $ \partial D $ is not smooth at $ p^0 $. We study $ F_D^k(z) $ as $ z \rightarrow p^0 $. It is evident from 
Definition \ref{ps_defn} that exactly two 
of the hypersurfaces $ \{ z \in \mathbb{C}^2: \rho_j(z) =  0\} $ (where $ j =1, \ldots, k $) intersect at the point $ p^0 $. Without loss of generality, we may assume that
\[
 \rho_1(p^0) = \rho_2(p^0) = 0.
\]
Let $ p^j $ be a sequence of points in $ D $ converging to $ p^0 $. Denote by
\begin{alignat*}{3}
 \lambda_j = \mbox{dist}\big( p^j, \{ \rho_1 = 0 \} \big), \\ 
 \mu_j = \mbox{dist}\big( p^j, \{ \rho_2 = 0 \} \big) 
\end{alignat*}
for each $ j $. Note that both $ \lambda_j $ and $ \mu_j $ tend to zero as $ j \rightarrow \infty $. 

\medskip
Following \cite{Kim-Yu}, there are three cases to be considered:
\begin{enumerate}
 \item[(I)] The sequence $ p^j $ is of radial type, i.e., there is a positive constant $ C $ (independent of $ j $) such that 
 $ 1/C \leq  \mu_j^{-1} {\lambda_j} \leq C $ for all $ j $.
\item[(II)] The sequence $ p^j $ is of $ q$-tangential type, i.e., either $ \lim_{j \rightarrow \infty} \mu_j^{-1} \sqrt{\lambda_j} = 0 $ or
$  \lim_{j \rightarrow \infty} \lambda_j^{-1} \sqrt{\mu_j} = 0 $.
\item[(III)] The sequence $ p^j $ is of mixed type, i.e., it is neither radial type nor $ q$-tangential type. Here, there are further two cases:
\begin{enumerate}
 \item [(a)] $ \lim_{j \rightarrow \infty} \frac{\lambda_j}{\mu_j} = 0 $ and $ \lim_{j \rightarrow \infty} \frac{\sqrt{\lambda_j}}{\mu_j} = m > 0 $,
 \item [(b)] $ \lim_{j \rightarrow \infty} \frac{\lambda_j}{\mu_j} = 0 $ and $ \lim_{j \rightarrow \infty} \frac{\sqrt{\lambda_j}}{\mu_j} = \infty $.
\end{enumerate}
\end{enumerate}

\begin{theorem} \label{ps_thm}
Let $ D $ be a strongly pseudoconvex polyhedral bounded domain in $ \mathbb{C}^2 $ with piecewise smooth boundary. 
Let $ p^0 \in \partial D $ be a singular boundary point and $ p^j $ be a sequence of points in $ D $ converging to $ p^0 $.
\begin{enumerate}
 \item [(i)] If the sequence $ \{p^j\} $ is of \textit{radial type}, then $  F^{k}_{D}(p^j) \rightarrow F^k_{\Delta \times \Delta} \big( (0,0)\big)= 1 $.
 \item [(ii)] If the sequence $ \{p^j\} $ is of \textit{q-tangential type}, then $  F^{k}_{D}(p^j) \rightarrow F^k_{\mathbb{B}^2} \big( (0,0)\big) =1 $.
 \item [(iii)] If the sequence $ \{p^j\} $ is of \textit{mixed type}, then 
 \[ F^{k}_{D}(p^j) \rightarrow F^k_{D_{1, \infty}} \big( (0,0)\big) \]
  in case (III)(a) and
 \[
  F^{k}_{D}(p^j) \rightarrow F^k_{ \Delta \times \Delta } \big( (0,0)\big) = 1 
  \]
  in case (III)(b),
\end{enumerate} 
where $ D_{1, \infty} $ is the model domain defined by 
\begin{equation*} 
D_{1, \infty} = \left\{ (z_1, z_2) \in \mathbb{C}^2: \Im z_1 + 1 > \frac{ Q_1(z_2)}{m^2}, \Im z_2 > -1 \right\}, 
\end{equation*}
and $ Q_1 $ is a strictly subharmonic polynomial of degree $ 2 $. 
\end{theorem}

It should be noted that if $ p^0 \in \partial D $ is a smooth boundary point, then the proof of Theorem 1.1 of \cite{BBMV} implies that 
$ F_D^k (z) \rightarrow 1 $ as $ z \rightarrow p^0 $.

\medskip

We adapt the scaling method from \cite{Kim-Yu} to understand $ F_D^k(p^j) $ in each of the above cases. To begin with, apply 
a complex linear change of coordinates $ A $ so that $ A(p^0) = (0,0) $ and the gradient vector to the hypersurface 
$ A \big(\{ \rho_1=0\} \big) $ and $ A \big(\{ \rho_2=0\} \big) $ at the origin is parallel to the $\Im z_1$ and $ \Im z_2$ axis respectively. Write 
$ A(p^j) = \tilde{p}^j $ for each $ j $.

\medskip

\noindent \textit{Case (I):} The smoothness of $ \rho_1 $ and $ \rho_2 $ implies that for each $ j $, there is a unique point
$ s^j $ on $ A\big( \{ \rho_1 = 0 \}\big) $ and $ t^j $ on $ A\big(\{ \rho_2 = 0 \}\big) $ such that 
\begin{alignat*}{3}
 \mbox{dist}\Big( \tilde{p}^j, A\big(\{ \rho_1 = 0 \}\big) \Big) = | \tilde{p}^j - s^j|,\\
\mbox{dist} \Big( \tilde{p}^j, A\big(\{ \rho_2 = 0 \}\big) \Big) = |\tilde{p}^j - t^j|.
\end{alignat*} 
There exists a sequence $ \{B^j \} $ of affine automorphisms of $ \mathbb{C}^2 $ such that $ B^j (\tilde{p}^j) = (0,0) $ for each $ j $ and the domains 
$ B^j \circ A (U \cap D) $ (for a sufficiently small neighbourhood $ U $ of $ p^0 $) are defined by 
\begin{multline*}
\Big\{ (z_1, z_2): \Im \big( z_1 - s^j_1 \big) > Q_1(z_2, \overline{z}_2) + o\big( |z_1 - s^j_1| + |z_2|^2\big),\\  \Im \big( z_2 - t^j_2 \big) 
> Q_2(z_1, \overline{z}_1) + o\big( |z_2 - t^j_2| + |z_1|^2\big) \Big\}, 
\end{multline*}
where $ Q_1 $ and $ Q_2 $ are real-valued quadratic polynomials. 

\medskip

Define the dilations 
\begin{alignat*}{3}
 L^j(z_1, z_2) = \left( \frac{z_1}{\lambda_j}, \frac{z_2}{\mu_j}\right),
\end{alignat*}
and the dilated domains $ D^j = L^j \circ B^j \circ A (U \cap D) $. Note that $ L^j \circ B^j \circ A (p^j) = (0,0) $ for all $ j $. 
Among other things, the following two claims were proved in \cite{Kim-Yu}. First, that $ D^j $ converges to 
\begin{alignat*}{3}
 D_{\infty} = \big\{ (z_1,z_2) \in \mathbb{C}^2: \Im z_1 > -c_1, \Im z_2 > -c_2 \big\}, 
\end{alignat*}
where $ c_1 $ and $ c_2 $ are positive constants. Secondly, for all $ j $ large, the scaled domains $ D^j $ are contained in $ D_0 $, where
\begin{equation*}
 D_0 = \big\{ (z_1,z_2) \in \mathbb{C}^2: \Im z_1 > -c_1 -r, \Im z_2 > -c_2 -r \big\},
\end{equation*}
and $ r > 0 $ is fixed. It should be noted that there is a biholomorphism from the limit domain $ D_{\infty} $ onto the unit bidisc $ \Delta \times \Delta $
that preserves the origin. 

\medskip 

\noindent \textit{Case (II):} Assume that the sequence $ p^j $ is of $ q$-tangential type to $ \{ \rho_1 = 0 \} $, 
i.e., $ \lim_{j \rightarrow \infty} \mu_j^{-1} \sqrt{\lambda_j} = 0 $. 

\medskip

For a sufficiently small neighbourhood $ U $ of $ p^0 $, we may assume that $ p^j $ are in $ U $ for all $ j $.  The domain $ A(U \cap D ) $ is given by
\begin{multline*}
 \Big\{ (z_1,z_2): \Im \big( z_1 - \tilde{p}^j_1 \big) + \lambda_j > Q_1 \big( z_2 - \tilde{p}^j_2 \big) + 
 o \big( \big| z_1 - \tilde{p}^j_1 \big| + \big| z_2 - \tilde{p}^j_2 \big|^2 \big), \\
 \Im \big( z_2 - \tilde{p}^j_2 \big) + \mu_j > Q_2 \big( z_1 - \tilde{p}^j_1 \big) + 
 o \big( \big| z_2 - \tilde{p}^j_2 \big| + \big| z_1 - \tilde{p}^j_1 \big|^2 \big)
 \Big\},
\end{multline*}
where $ Q_1 $ are $ Q_2 $ are strictly subharmonic quadratic polynomials. Let $ L^j : \mathbb{C}^2 \rightarrow \mathbb{C}^2 $ be the dilations given by
\begin{equation*}
 L^j (z_1, z_2) = \left( \frac{z_1 - \tilde{p}^j_1}{ \lambda_j}, \frac{z_1 - \tilde{p}^j_2}{ \sqrt{\lambda_j}} \right).
\end{equation*}
It follows that $ L^j \circ A (p^j) = (0, 0) $ and the scaled domains $ D^j = L^j \circ A (U \cap D) $ converge to 
\begin{equation*}
 D_{\infty} = \big\{ (z_1, z_2) \in \mathbb{C}^2 : \Im z_1 + 1 > Q_1 (z_2) \big\}, 
\end{equation*}
which is biholomorphically equivalent to $ \mathbb{B}^2 $.

\medskip 

\noindent \textit{Case (III):} Here, the sequence $ p^j $ is of mixed type. Consider the dilations  
\begin{alignat*}{3}
 L^j (z_1, z_2) = \left( \frac{z_1 - \tilde{p}^j_1}{ \lambda_j}, \frac{z_1 - \tilde{p}^j_2}{ {\mu_j}} \right)
\end{alignat*}
and note that $ L^j \circ A (p^j) = (0, 0) $. It follows that the dilated domains $ D^j = L^j \circ A (U \cap D) $ converge to 
\begin{alignat*}{3}
 D_{\infty} = \left\{ (z_1, z_2) \in \mathbb{C}^2: \Im z_1 + 1 > \lim_{j \rightarrow} \frac{ \mu_j^2}{\lambda_j} Q_1(z_2), \Im z_2 > -1 \right\}.
\end{alignat*}
More specifically, the limit domain turns out to be 
\begin{equation} \label{E1}
D_{1, \infty} = \left\{ (z_1, z_2) \in \mathbb{C}^2: \Im z_1 + 1 > \frac{ Q_1(z_2)}{m^2}, \Im z_2 > -1 \right\}.
\end{equation}
in case III(a), and 
\begin{alignat*}{3} \label{E2}
D_{2,\infty} = \left\{ (z_1, z_2) \in \mathbb{C}^2: \Im z_1 > -1, \Im z_2 > -1 \right\} 
\end{alignat*}
in case III(b). 

Note that the limiting domain $ D_{1, \infty} $ is a Siegel domain of second kind (refer \cite{Piatetski-Shapiro} for more details) and 
hence complete Kobayashi hyperbolic. Evidently, $ D_{1, \infty} $ can be written as the intersection of 
an open ball with a half space in $ \mathbb{C}^2 $. Moreover, $ D_{1,\infty} $ is an unbounded convex domain. Furthermore, according to \cite{Piatetski-Shapiro}, $ D_{1,\infty} $ is 
biholomorphic to a bounded domain in $ \mathbb{C}^2 $. In particular, the Bergman kernel $ K_{D_{1,\infty}} $ is non-vanishing along the 
diagonal. Also, note that the limit domain $ D_{2,\infty} $ is biholomorphic to the unit bidisc $ \Delta \times \Delta $ via a map that preserves the origin.

The stability of the infinitesimal Kobayashi metric under scaling can be proved using similar ideas as in Lemma 5.2 of \cite{MV}. The following two ingredients will be required
in the proof -- first, the limit domain $ D_{\infty} $ is complete Kobayashi hyperbolic and hence taut in each of the cases (I), (II), and (III). The 
next step is to consider the mappings $ f^j : \Delta \rightarrow D^j $ that almost realize $ k_{D^j} ( \cdot, \cdot) $ and establish that $ \{f^j\} $ is normal.
Recall that, in each of the three cases listed above, the scaled domains $ D^j $ are all contained in the taut domain $ 2 D_{\infty} $ for large $j $. 
Hence, it is possible to pass to a subsequence of $ \{f^j \} $ that converges to a
holomorphic mapping $ f : \Delta \rightarrow D_{\infty} $ uniformly on compact sets of $ \Delta $. It follows that the limit
map $ f $ provides a candidate in the definition of $ k_{D_{\infty}} (\cdot, \cdot) $.

\begin{lemma} \label{met-conv}

For $ (z,v) \in D_{\infty} \times \mathbb{C}^2 $, 
\begin{equation*}
  k_{D^j}(z,v) \rightarrow k_{D_{\infty}} (z,v).
\end{equation*}
Moreover, the convergence is uniform on compact sets of $ D_{\infty} \times \mathbb{C}^2 $.
\end{lemma}

The next step is a stability statement for the Kobayashi indicatrices of the scaled domains $ D^j $. 

\begin{lemma} \label{ind-conv} 
For $z$ in any compact subset $S$ of $D_{\infty} $,
\begin{enumerate}
\item [(i)] $ I^{k}_{D^j} (z) $ is uniformly compactly contained in $ \mathbb{C}^n $ for all $ j $ large, and
\item [(ii)] the indicatrices $ I^{k}_{D^j}(z) $ converge uniformly in the Hausdorff sense to $ I^{k}_{D_{\infty}}(z)$.
\end{enumerate}
Finally, for each $z \in D_{\infty}$, the functions $\la\big(I^{k}_{D^j}(z)\big)$ converge to $\la\big(I^{k}_{D_{\infty}}(z)\big)$.
\end{lemma}
For the proof, repeat the arguments provided earlier along with the following observation: the limit domain $ D_{\infty} $ is biholomorphically equivalent to a bounded domain in $ \mathbb{C}^2 $ in each of the cases (I),
(II) and (III), which implies that there is a uniform positive constant  $ C $ (depending only on $ S $) such that for $ z \in S $
\[
k_{D_\infty}(z,v) \geq C \vert v\vert
\]
for all $v \in \mathbb{C}^2$.

\textit{Proof of} Theorem \ref{ps_thm}: Observe that
\begin{alignat*}{3}
 F^{k}_{U \cap D}(p^j)=F^{k}_{D^j} \left( (0,0) \right) = K_{D^j} \big( (0,0) \big) \la\big(I^{k}_{D^j} (0,0)\big)
\end{alignat*}
for each $j $. To control the Bergman kernels $ K_{D^j} $ on the scaled domains, note first that the limit domain $ D_{\infty} $ is biholomorphic to a convex domain in each of the cases (I), 
(II) and (III)
which implies that
\begin{alignat}{3} \label{E4}
 K_{D^j} \big( (0,0) \big) \rightarrow K_{D_{\infty}} \big( (0,0) \big) 
\end{alignat}
by virtue of Lemma 2.1 of \cite{BBMV}. Moreover, applying Lemma \ref{ind-conv}, it follows that 
\begin{alignat}{3} \label{E3}
 F^{k}_{U \cap D}(p^j) \rightarrow F^k_{D_{\infty}} \big( (0,0)\big).
\end{alignat}
Finally, to conclude, note that the domain $ D $ as in Definition \ref{ps_defn} supports a local holomorphic peak function at each boundary point. It follows that $ F^k_D $ 
can be localized near $ p^0 \in \partial D $. This observation together with \eqref{E3} yields
\begin{alignat*}{3} 
 F^{k}_{D}(p^j) \rightarrow F^k_{D_{\infty}} \big( (0,0)\big),
\end{alignat*}
where $ D_{\infty} $ is the model domain at the point $ p^0 $. The result follows by recalling that the limit domain $ D_{\infty} $ is 
biholomorphic to $ \Delta \times \Delta $ in cases (I) and (IIIb) and to $ \mathbb{B}^2 $ in case (II).
\qed


\subsection*{Acknowledgements} The second named author was partially supported by the DST-INSPIRE grant IFA-13 MA-21.


\end{document}